\def\Thst{1.1}
\def\LeRi{1.2}
\def\NullD{1.3}

\def\Multt{2.3}

\def\Igz{3.1}
\def\RiGaz{3.2}
\def\tilAG{4.1}
\def\ActGn{4.2}

\def\Trafsp{4.4}
\def\VaR{4.5}
\catcode`\^^Z=9
\catcode`\^^M=10
\output={\if N\header\headline={\hfill}\fi
\plainoutput\global\let\header=Y}
\magnification\magstep1
\tolerance = 500
\hsize=14.4true cm
\vsize=22.5true cm
\parindent=6true mm\overfullrule=2pt
\newcount\kapnum \kapnum=0
\newcount\parnum \parnum=0
\newcount\procnum \procnum=0
\newcount\nicknum \nicknum=1
\font\ninett=cmtt9

\font\ninebf=cmbx9

\font\sixbf=cmbx6
\font\ninesl=cmsl9

\font\nineit=cmti9

\font\ninerm=cmr9

\font\sixrm=cmr6
\font\ninei=cmmi9
\font\eighti=cmmi8
\font\sixi=cmmi6
\skewchar\ninei='177 \skewchar\eighti='177 \skewchar\sixi='177
\font\ninesy=cmsy9
\font\eightsy=cmsy8
\font\sixsy=cmsy6
\skewchar\ninesy='60 \skewchar\eightsy='60 \skewchar\sixsy='60
\font\titelfont=cmr10 scaled 1440
\font\paragratit=cmbx10 scaled 1200

\font\name=cmcsc10
\font\emph=cmbxti10

\font\tenmsbm=msbm10
\font\sevenmsbm=msbm7
%

%
\font\got=eufm10
\font\Got=eufm7
\font\teneufm=eufm10
\font\seveneufm=eufm7
\font\fiveeufm=eufm5
\newfam\eufmfam
\textfont\eufmfam=\teneufm
\scriptfont\eufmfam=\seveneufm
\scriptscriptfont\eufmfam=\fiveeufm

\font\tenmsam=msam10
\font\sevenmsam=msam7
\font\fivemsam=msam5
\newfam\msamfam
\textfont\msamfam=\tenmsam
\scriptfont\msamfam=\sevenmsam
\scriptscriptfont\msamfam=\fivemsam
\font\tenmsbm=msbm10
\font\sevenmsbm=msbm7
\font\fivemsbm=msbm5
\newfam\msbmfam
\textfont\msbmfam=\tenmsbm
\scriptfont\msbmfam=\sevenmsbm
\scriptscriptfont\msbmfam=\fivemsbm
\def\Bbb#1{{\fam\msbmfam\relax#1}}
\def\cz{{\kern0.4pt\Bbb C\kern0.7pt}
}
\def\ez{{\kern0.4pt\Bbb E\kern0.7pt}
}
\def\fz{{\kern0.4pt\Bbb F\kern0.3pt}}
\def\gz{{\kern0.4pt\Bbb Z\kern0.7pt}}
\def\hz{{\kern0.4pt\Bbb H\kern0.7pt}
}
\def\kz{{\kern0.4pt\Bbb K\kern0.7pt}
}
\def\nz{{\kern0.4pt\Bbb N\kern0.7pt}
}
\def\oz{{\kern0.4pt\Bbb O\kern0.7pt}
}
\def\rz{{\kern0.4pt\Bbb R\kern0.7pt}
}
\def\sz{{\kern0.4pt\Bbb S\kern0.7pt}
}
\def\pz{{\kern0.4pt\Bbb P\kern0.7pt}
}
\def\qz{{\kern0.4pt\Bbb Q\kern0.7pt}
}
\newskip\ttglue
\def\ninepoint{\def\rm{\fam0\ninerm}%
  \textfont0=\ninerm \scriptfont0=\sixrm \scriptscriptfont0=\fiverm
  \textfont1=\ninei \scriptfont1=\sixi \scriptscriptfont1=\fivei
  \textfont2=\ninesy \scriptfont2=\sixsy \scriptscriptfont2=\fivesy
  \textfont3=\tenex \scriptfont3=\tenex \scriptscriptfont3=\tenex
  \def\it{\fam\itfam\nineit}%
  \textfont\itfam=\nineit
  \def\sl{\fam\slfam\ninesl}%
  \textfont\slfam=\ninesl
  \def\bf{\fam\bffam\ninebf}%
  \textfont\bffam=\ninebf \scriptfont\bffam=\sixbf
   \scriptscriptfont\bffam=\fivebf
  \def\tt{\fam\ttfam\ninett}%
  \textfont\ttfam=\ninett
  \tt \ttglue=.5em plus.25em minus.15em
  \normalbaselineskip=11pt
  \font\name=cmcsc9
  \let\sc=\sevenrm
  \let\big=\ninebig
  \setbox\strutbox=\hbox{\vrule height8pt depth3pt width0pt}%
  \normalbaselines\rm
  \def\sl{\it}}

\headline={\ifodd\pageno\rightheadline\else\leftheadline\fi}
\def\rightheadline{\ninepoint Paragraphen"uberschrift\hfill\folio}
\def\leftheadline{\ninepoint\folio\hfill Chapter"uberschrift}
\let\header=Y
\def\titel#1{\need 9cm \vskip 2truecm
\parnum=0\global\advance \kapnum by 1
{\baselineskip=16pt\lineskip=16pt\rightskip0pt
plus4em\spaceskip.3333em\xspaceskip.5em\pretolerance=10000\noindent
\titelfont Chapter \uppercase\expandafter{\romannumeral\kapnum}.
#1\vskip2true cm}\def\leftheadline{\ninepoint
\folio\hfill Chapter \uppercase\expandafter{\romannumeral\kapnum}.
#1}\let\header=N
}
\def\Titel#1{\need 9cm \vskip 2truecm
\global\advance \kapnum by 1
{\baselineskip=16pt\lineskip=16pt\rightskip0pt
plus4em\spaceskip.3333em\xspaceskip.5em\pretolerance=10000\noindent
\titelfont\uppercase\expandafter{\romannumeral\kapnum}.
#1\vskip2true cm}\def\leftheadline{\ninepoint
\folio\hfill\uppercase\expandafter{\romannumeral\kapnum}.
#1}\let\header=N
}
\def\need#1cm {\par\dimen0=\pagetotal\ifdim\dimen0<\vsize
\global\advance\dimen0by#1 true cm
\ifdim\dimen0>\vsize\vfil\eject\noindent\fi\fi}
\def\neupara#1{\par\penalty-2000
\procnum=0\global\advance\parnum by 1
\vskip1cm\noindent{\paragratit \the\parnum. #1}%
\def\rightheadline{\ninepoint\S\the\parnum.\ #1\hfill \folio}%
\vskip 8mm\noindent}
\def\Proclaim #1 #2\finishproclaim {\bigbreak\noindent
{\bf#1\unskip{}. }{\it#2}\medbreak\noindent}
%
\gdef\proclaim #1 #2 #3\finishproclaim {\bigbreak\noindent%
\global\advance\procnum by 1
{%
{\relax\ifodd \nicknum
\hbox to 0pt{\vrule depth 0pt height0pt width\hsize
   \quad \ninett#3\hss}\else {}\fi}%
\bf\the\parnum.\the\procnum\ #1\unskip{}. }
{\it#2}
\medbreak\noindent}
\newcount\stunde \newcount\minute \newcount\hilfsvar
\def\uhrzeit{
    \stunde=\the\time \divide \stunde by 60
    \minute=\the\time
    \hilfsvar=\stunde \multiply \hilfsvar by 60
    \advance \minute by -\hilfsvar
    \ifnum\the\stunde<10
    \ifnum\the\minute<10
    0\the\stunde:0\the\minute~Uhr
    \else
    0\the\stunde:\the\minute~Uhr
    \fi
    \else
    \ifnum\the\minute<10
    \the\stunde:0\the\minute~Uhr
    \else
    \the\stunde:\the\minute~Uhr
    \fi
    \fi
    }

 \def\calT{{\cal T}}

\def\gotn{\hbox{\got n}}

\def\Gotn{\hbox{\Got n}}

\def\GL{\mathop{\rm GL}\nolimits}

\def\kernel{\mathop{\rm kernel}\nolimits}

\def\mod{\mathop{\rm mod}\nolimits}

\def\SL{\mathop{\rm SL}\nolimits}

\def\Sp{\mathop{\rm Sp}\nolimits}

\def\boxit#1{
  \vbox{\hrule\hbox{\vrule\kern6pt
  \vbox{\kern8pt#1\kern8pt}\kern6pt\vrule}\hrule}}
\def\Boxit#1{
  \vbox{\hrule\hbox{\vrule\kern2pt
  \vbox{\kern2pt#1\kern2pt}\kern2pt\vrule}\hrule}}

\def\zwischen#1{\bigbreak\noindent{\bf#1\medbreak\noindent}}

\def\smallni{\smallskip\noindent }
\def\medni{\medskip\noindent }
\def\bigni{\bigskip\noindent }
\def\Isom{\mathop{\;{\buildrel \sim\over\longrightarrow }\;}}
\def\lo{\longrightarrow}

\def\loma{\longmapsto}
\def\betr#1{\vert#1\vert}

\def\imag{{\rm i}}
\def\pii{\pi {\rm i}}

\def\set#1{\bigl\{\,#1\,\bigr\}}

\def\square{\hbox{\hbox to 0pt{$\sqcup$\hss}\hbox{$\sqcap$}}}
\def\qed{\ifmmode\square\else{\unskip\nobreak\hfil
\penalty50\hskip3em\null\nobreak\hfil\square
\parfillskip=0pt\finalhyphendemerits=0\endgraf}\fi}
\def\pn{\the\parnum.\the\procnum}
\def\downmapsto{{\buildrel
        {\vbox{\hbox{\hskip.2pt$\scriptstyle-$}}}
        \over{\raise7pt\vbox{\vskip-4pt\hbox{$\textstyle\downarrow$}}}}}
\def\Hilb{\hbox{\rm Hilb}}
\nopagenumbers
\immediate\newwrite\num
\nicknum=0  
\let\header=N

\immediate\openout\num=calabi1.num
\immediate\newwrite\num\immediate\openout\num=calabi1.num
\def\RAND#1{\vskip0pt\hbox to 0mm{\hss\vtop to 0pt{%
  \raggedright\ninepoint\parindent=0pt%
  \baselineskip=1pt\hsize=2cm #1\vss}}\noindent}
\noindent
\centerline{\titelfont Some Siegel threefolds with a Calabi-Yau model}%
\def\leftheadline{\ninepoint\folio\hfill
Some Siegel threefolds with a Calabi-Yau model}%
\def\rightheadline{\ninepoint Introduction\hfill \folio}%
\headline={\ifodd\pageno\rightheadline\else\leftheadline\fi}

\vskip 1.5cm
\leftline{\it \hbox to 6cm{Eberhard Freitag\hss}
Riccardo Salvati
Manni  }
  \leftline {\it  \hbox to 6cm{Mathematisches Institut\hss}
Dipartimento di Matematica, }
\leftline {\it  \hbox to 6cm{Im Neuenheimer Feld 288\hss}
Piazzale Aldo Moro, 2}
\leftline {\it  \hbox to 6cm{D69120 Heidelberg\hss}
 I-00185 Roma, Italy. }
\leftline {\tt \hbox to 6cm{freitag@mathi.uni-heidelberg.de\hss}
salvati@mat.uniroma1.it}
\vskip1cm
\centerline{\paragratit \rm  2009}%
\vskip5mm\noindent%
\let\header=N%
\def\imag{{\rm i}}%
{\paragratit Introduction}%
\medni
In the following we describe some examples of Calabi-Yau manifolds that
arise as desingularizations of certain Siegel threefolds.
Here by a Calabi-Yau manifold we understand a smooth complex projective
variety which admits a holomorphic differential form of degree three
without zeros and such that the first Betti number is zero.
This differential form is
unique up to a constant factor, and we call it the Calabi-Yau form.
Our  interest in this subject is influenced by work of Gritsenko and
many discussions with him.
The first Siegel modular variety with a Calabi-Yau model and the essentially only one
up to now has been discovered by Barth and Nieto. They showed that the
``Nieto quintic''
$\{x\in P^5(\cz),\ \sigma_1(x)=\sigma_5(x)=0\}$, where $\sigma_i$ denote the elementary
symmetric polynomials, has a Calabi-Yau model and they derived that the Siegel
modular variety $A_{1,3}(2)$ of polarization type $(1,3)$ and a certain level two structure has a
Calabi-Yau model.
Since the Jacobian of a symplectic substitution is $\det(CZ+D)^{-3}$,
the Calabi-Yau three-form produces a modular  form of weight three
and this must be a cusp form, since it survives on a non-singular model as a
holomorphic differential form ([Fr], III.2.6).
In the paper [GH] Gritsenko and Hulek gave a direct construction of this modular form
and obtained a new proof for the fact that $A_{1,3}(2)$ has a Calabi-Yau model.
We also refer to [GHSS] for further investigations.
Besides this example and some small extensions of this group with the same three-form
no other examples of Siegel threefolds with Calabi-Yau model seem to be known.
Gritsenko raised the problem of determing
all Siegel threefolds which admit a Calabi-Yau model.
As we mentioned already, such a threefold will produce a certain cusp form
of weight three for the considered modular group $\Gamma$.
This cusp form has very restrictive properties. Since the induced
differential form should have no zero at least at the regular locus of the
quotient $\hz_2/\Gamma$, all zeros of the form must be contained in the ramification
of $\hz_2\to\hz_2/\Gamma$. Gritsenko gave examples of such modular forms:
we refer to the paper [GC] which contains some systematic study of them. One example 
that
Gritsenko and Cl\'ery describe, is the form $\nabla_3$, which
is a cusp form
of weight three for the Hecke group
$\Gamma_{2,0}[2]$ with respect to a certain quadratic character $\chi$.
Hence a subgroup of index
two is a candidate for producing a Calabi-Yau manifold.
We will prove that this is the case. We will show more: for any group
between $\Gamma_2[2]$ and  $\Gamma_{2,0}[2]$ there exists a subgroup of index
two (the kernel of $\chi$) which produces a Calabi-Yau manifold.
\smallskip
The modular form $\nabla_3$ will come up in a completely different
manner. It is simply  the product $T$ of 6 (of the 10
classical) theta constants with suitable
properties. In this form it has already been described  in [GS]
and these expressions occur also in [GN].
This approach has the advantage that we can easily  describe the action of the
full modular group, which is necessary, since we need information about this
form at all boundary components.
Another advantage of this description is that we can
use the work  of Igusa about the structure
of the ring of modular forms with respect to his group $\Gamma_2[4,8]$
and of some groups containing this group.
Igusa used the ten theta constants of the first kind.
If one is concerned with groups in the region of the principal
congruence subgroup of level two, there are advantages in using the theta constants
of the second kind. We use the very nice approach  given by
Runge [Ru1], [Ru2].
\smallskip
We also make use of
Igusa's method of desingularization of the Siegel threefold
with respect to the principal congruence subgroup
of level $q>2$ (we need $q=4$).
\smallskip
Using Igusa's results or Runge's approach, it is easy
to determine the
rings
of modular forms for the groups in question, and in this way one can produce
equations for the Siegel threefolds.
The main example is the subgroup of index two of $\Gamma_2[2]$.
\smallskip
In this introduction we only describe the equations of this
Siegel threefold  in a purely algebraic
way.
\Proclaim
{Theorem}
{Let $X$ be the subvariety of $P^5(\cz)$ given by the intersection
of the quartic
$$y_5^4=y_0y_1y_2( y_0+y_1+y_2+ y_3+y_4)$$
and the quadric
$$2y_5^2= y_0y_1+y_0y_2+ y_1y_2- y_3y_4.$$
This is a normal projective variety
of dimension three. There exists a desingularization
$\tilde X\to X$ which is a Calabi-Yau manifold.
}
\finishproclaim
The variety $X$ together
with the Calabi-Yau form have some symmetries.
They are easier to describe in another coordinate
system (see \VaR\ for the explicit description).
In this coordinate system we also will give an
explicit algebraic expression for the Calabi-Yau form.
We will see:
\Proclaim
{Theorem}
{There is a subgroup $G\subset PGL(5,\cz)$,
isomorphic to the semidirect product
$S_3\cdot (\gz/2\gz)^2$, which leaves $X$ and the form $\omega$
invariant.
For each subgroup
$H$ of $G$
the quotient $X/H$ admits a desingularization
which is a Calabi-Yau manifold.}
\finishproclaim
It is not difficult to determine the singular locus of $X$:
\Proclaim
{Proposition}
{The singular locus of $X$ is the union of $15$ smooth curves.
It consists of two $G$-orbits. One orbit consists of three quadrics, the other
of
$12$
lines. Representatives are given by the ideals
$$(y_0+y_4,y_1+y_4,y_3-y_4,y_2y_4+y_5^2),\quad
(y_0,y_2,y_3,y_5).$$
}
\finishproclaim
The main problem is to find a good resolution of the singularities of $X$.
It might be possible to do this by hand or with the help of
a computer.
We will find it by interpreting $X$ as a Siegel modular variety.\bigskip

{\paragratit Acknowledgements}%

\medni
We would like to thank T. Bridgeland, A.Corti, S.Cynk  and  C.Meyer  for useful discussions

\neupara{The Siegel modular group of genus two}%
Recall that the real symplectic group
$$\Sp(n,\rz)=\set{M\in\GL(2n,\rz);
\quad {^tM}IM=I}\qquad \left(I=\pmatrix{0&-E\cr E&0}\right)$$
acts on the generalized half plane
$$\hz_n:=\set{Z=X+\imag Y;\quad Z={^tZ},\  Y>0\ \hbox{(positive definite)}}$$
by
$$MZ=(AZ+B)(CZ+D)^{-1},\qquad M=\pmatrix{A&B\cr C&D}.$$
Let $\Gamma_n:=\Sp(n,\gz)$ be  the Siegel modular group.
The principal congruence subgroup of level $l$ iy
$$\Gamma_n[l]:=\kernel (\Sp(n,\gz)\lo\Sp(n,\gz/l\gz))$$
and Igusa's subgroup is
$$\Gamma_n[l,2l]:=\set{M\in\Gamma_n[l];\quad
A{^tB}/l\ \hbox{and}\ C{^tD}/l\ \hbox{have even diagonal}}.$$
For even $l$, Igusa's subgroup  is a normal subgroup of $\Gamma_n$. Other important subgroups
are
$$\Gamma_{n,0}[l]=\{M\in\Gamma_n;\quad C\equiv 0\;\mod\; l\}.$$
\eject
\zwischen{Theta characteristics in genus two}%
By definition, a theta characteristic is an element $m={a\choose b}$ from
$(\gz/2\gz)^{2n}$. Here $a,b\in(\gz/2\gz)^n$. The characteristic is called even
if ${^ta}b=0$.
The group $\Sp(n,\gz/2\gz)$ acts on the set of characteristics by
$$M\{m\}:={^tM}\strut^{-1}m+\pmatrix{(A{^tB})_0\cr (C{^tD})_0}.$$
Here $S_0$ denotes the column built of the diagonal of a square matrix $S$.
It is well-known that $\Sp(n,\gz/2\gz)$ acts transitively on the subsets of
even and odd characteristics.
Recall that for any characteristic the theta function
$$\vartheta[m]=\sum_{g\in\gz^n}e^{\pii (Z[g+a/2]+{^tb}(g+a/2))}$$
can be defined.
Here we use the identification of $\gz/2\gz$ with the
subset $\{0,1\}\subset\gz$. It vanishes if and only if $m$ is odd.
Recall also that the formula
$$\vartheta[M\{m\}](MZ)=v(M,m)\sqrt{\det (CZ+D)}\vartheta[m](Z)$$
holds for $M\in\Gamma_n$, where $v(M,m)$ is a rather delicate
$8^{\hbox{\sevenrm th}}$ root of unity which depends on the choice of the
square root.
\smallskip
We are interested in the case $n=2$. In this case there are ten even characteristics.
We will write the coordinates of $\hz_2$ as
$$Z=\pmatrix{z_0&z_1\cr z_1&z_2}.$$
A set $\{m_1,m_2,m_3,m_4\}$ of four
pairwise different even
characteristics is called {\it syzygetic\/} if the sum of any three of them is even.
There are 15 syzygetic (unordered) quadruples and the group $\Sp(2,\gz)$
acts transitively on them. We call
$$\left\{\matrix{ 0&0&0&0\cr 0&0&0&0\cr 0&1&0&1\cr 0&0&1&1}\right\}$$
the standard syzygetic quadruple. We are also interested in the 15 complementary
sextuples of even characteristics $\{n_1,\dots,n_6\}$.
From a detailed study of the multipliers $v(M,m)$ one can deduce:
\proclaim
{Lemma}
{Let  $\gotn=\{n_1,\dots,n_6\}$ be a sextuple of even
characteristics (in the case $n=2$)
complementary to a syzygetic quadruple. Then the product
$$T=T_{\Gotn}:=\prod_{\nu=1}^6\vartheta[n_\nu](Z)$$
is a cusp form of weight\/ $3$ for a group 
conjugate to  $\Gamma_{2,0}[2]$ and with respect
to a quadratic character $\chi_{\Gotn}$ on this group.
The kernel of this character
contains $\Gamma_2[4]$.
In the case of the standard tuple the group is $\Gamma_{2,0}[2]$
and the kernel
$$\Gamma_{2,0}[2]_{\Gotn}:=\bigl\{\;M\in\Gamma_{2,0}[2],\quad \chi_{\Gotn}(M)=1\;\bigr\}$$
is a subgroup of index two.
}
Thst%
\finishproclaim
We describe the character of $T$ for the standard syzygetic quadruple.
For this we introduce
$$\Theta:= \vartheta\Bigl[{00\atop01}\Bigr] \vartheta\Bigl[{00\atop00}\Bigr]
\vartheta\Bigl[{00\atop10}\Bigr] \vartheta\Bigl[{00\atop11}\Bigr].$$
Since $T\cdot\Theta=\chi_5$ is a modular form with respect to the full modular
group, this form also is a modular form for $\Gamma_{2,0}[2]$. The characters
of $T$ and $\Theta$ differ by the character of $\chi_5$.
Using the well-known isomorphism $\Gamma_2/\Gamma_2[2]\cong S_6$,
the character of $\chi_5$ corresponds to the sign character of $S_6$.
Hence it is sufficient
to describe the character of $\Theta$.
\proclaim
{Lemma}
{The character of $\Theta$ for $M\in\Gamma_{2,0}[2]$ is given by
$$(-1)^{(\alpha+\beta+ \gamma)/2},$$
where
$$C\,{\strut^tD} = \pmatrix{\alpha&\beta\cr \beta& \gamma}.$$
The character of $T$ is the product of this character and the only nontrivial
character of the full modular group.
}
LeRi%
\finishproclaim
The proof can be taken from [SM] (use Lemma 4).
\smallskip
We have to recall the location of the zeros of the theta functions.
In the case ${^tm}=(1,1,1,1)$ the theta function $\vartheta[m]$ has a zero
of first order along
the diagonal $z_1=0$ and every zero component is equivalent to the diagonal
mod $\Gamma_2[1,2]$.
We consider the differential form
$$\omega=\omega_{\Gotn}=T_{\Gotn}\;dz_0\wedge dz_1\wedge dz_2.$$
It is invariant under the group $\Gamma_{2,0}[2]_{\Gotn}$.
\proclaim
{Proposition}
{let $F\subset\hz_2$ be an irreducible component of the zero
locus of $\omega_{\Gotn}$. Then there exists an element
$M\in\Gamma_2$
whose fixed point set is $F$. It has the property $M^2=E$.
This element is uniquely determined
up to the sign and is actually contained already in
$$\Gamma_2[2,4]_{\Gotn}=\Gamma_2[2,4]\cap\Gamma_{2,0}[2]_{\Gotn}.$$
}
NullD%
\finishproclaim
{\it Proof.\/} Taking a conjugate group we may assume that $F$ is the diagonal
$z_1=0$. Then in the sextuple the characteristic  ${^tm}=(1,1,1,1)$ must
occur.
The only non trivial modular substitution which fixes the diagonal
is $z_1\mapsto-z_1$. The theta series for ${^tm}=(1,1,1,1)$ changes its sign
under this substitution. The others are invariant.
This shows
$$T_{\Gotn}\pmatrix{z_0&-z_1\cr -z_1&z_2}=-T_{\Gotn}\pmatrix{z_0&z_1\cr z_1&z_2}.$$
Since this is the transformation law for a modular form of odd weight, the
substitution $z_1\mapsto-z_1$ is in the kernel of $\chi_{\Gotn}$.
It is also in $\Gamma_2[2,4]$.\qed
\smallskip
We have to consider the group
$$\Gamma=\Gamma_{\Gotn}:=\bigl\{M\in\Gamma_2[2];\quad
\chi_{\Gotn}(M)=1\bigr\}.$$
This is a subgroup of index two of $\Gamma_2[2]$.
Using \LeRi\ and some computation one can see:
\proclaim
{Lemma}
{
For the standard syzygetic quadruple, the group
$\Gamma$ is defined inside $\Gamma_2[2]$  by the condition that
the symmetric matrix
$C\,{\strut^tD} = {\alpha\;\beta\choose \beta\; \gamma}$
has the property
$\alpha+\beta+ \gamma\equiv 0$ $\mod4$. Moreover 
$\Gamma=\Gamma_{\Gotn}\subset\Gamma_2[2]$ is generated by
\vskip1mm
\item{\rm 1)} The group $\Gamma_2[4]$,
\item{\rm 2)} The elements of $\Gamma_{2,0}[2]_{\Gotn}$ that are
conjugate inside $\Gamma_2$ to the diagonal matrix with diagonal
$(1,-1,1,-1)$.
\item{\rm 3)} All elements of $\Gamma_{2,0}[2]_{\Gotn}$ that are
conjugate inside $\Gamma_2$ to a translation matrix ${E\,S\choose0\,E}$ of
$\Gamma_2[2]$.
\vskip0pt
}
DGa%
\finishproclaim
\neupara{Igusa's desingularization}%
We consider the principal congruence subgroup $\Gamma_2[l]$ of level $l\ge 3$
and denote by
$$X=X(l):=\overline{\hz_2/\Gamma_2[l]}$$
the Satake compactification and by $\tilde X=\tilde X(l)$ the monoidal
transform along the Satake boundary. Igusa proved that $\tilde X$ is smooth.
The theory of Igusa is very difficult but fortunately we can formulate
in a very  simple way what we need from it.
\smallskip
Igusa used so-called normal
coordinates
$$q_0=e^{2\pii (z_0+z_1)/l},\quad q_2=e^{2\pii (z_2+z_1)/l},
\quad q_1=e^{-2\pii z_1/l}.$$
We consider them as a holomorphic map
$$\hz_2\lo \cz^3,\quad \pmatrix{z_0&z_1\cr z_1& z_2}\loma (q_0,q_1,q_2).$$
The image is an open subset $D\subset \cz^3$.
If we denote by $\calT=\calT(l)$ the group of translations $Z\mapsto Z+lS$,
$S$ integral, we get a biholomorphic map
$$\hz_2/\calT\Isom D.$$
We can consider its inverse map and compose it with the projection onto
$\hz_2/\Gamma_2[l]$ and the inclusion into $X$ to get a holomorphic
map
$$D\lo X.$$
A domain $U\subset\cz^n$ is called a {\it Reinhardt domain\/} if
$$(z_1,\dots,z_n)\in U\Longrightarrow (\zeta_1 z_1,\dots,\zeta_n z_n)\in U
\quad\hbox{for}\ \betr{\zeta_\nu}=1.$$
It is called a {\it complete Reinhardt domain,\/} if this is true for all $\zeta$ with
$\betr{\zeta_\nu}\le 1$. Each Reinhardt domain can be completed to a complete
Reinhardt domain in an obvious way. Any holomorphic function on a Reinhardt domain
can be expanded  in a Laurent series in the whole domain. This Laurent series is
a power series if and only if the function extends as a holomorphic function to the
completed Reinhardt domain.
\proclaim
{Lemma}
{The domain $D$ is a Reinhardt domain $\cz^3$. Its completion is
$$\tilde D:=D\cup\{q\in\cz^3;\ q_0q_1q_2=0\}.$$
The domain $D$ is dense in $\tilde D$.
}
RD%
\finishproclaim
All we need from Igusa's theory is:
\proclaim
{Theorem (Igusa)}
{The natural
map $D\to X$ extends to a locally biholomorphic map
$$\tilde D\lo\tilde X.$$
The group $\Gamma_2/\Gamma_2[l]$ acts on $X$ and hence on $\tilde X$. The translates of the
images of $\tilde D$ cover $\tilde X$.}
Igus%
\finishproclaim
Each holomorphic function on $D$ can be written as Laurent series in the variables
$q_\nu$.
Assume that it is the Fourier expansion of a modular form
$$\sum_T a(T)e^{2\pii\sigma(t_0z_0+2t_1z_1+t_2z_2)/l},\quad 
T=\pmatrix{t_0&t_1\cr t_1&t_2}.$$
Here $T$ runs over all  matrices such that $t_0,t_2,2t_1$
are integral. If we rewrite this in the normal coordinates, we get
$$\sum_T a(T) q_0^{t_0}q_1^{t_0+t_2-2t_1}q_2^{t_2}.$$
By  Koecher's principle we can have $a(T)\ne 0$ only for semi-positive $T$.
This means $t_0,t_2\ge 0$ and $t_0t_2\ge t_1^2$. This implies
$t_0+t_2-2t_1\ge 0$. Hence the Laurent series  actually is a power series.
Hence this function extends to a holomorphic
function on the whole $\tilde D$.
\smallskip
We can talk about the multiplicity of the zero along each of the $q_\nu=0$.
The theta functions $\vartheta[m]$ are periodic under $Z\mapsto 8S$,
$S$ integral. Hence we can use $l=8$ for the normal coordinates. The multiplicities
are easy to compute:
\proclaim
{Lemma}
{The multiplicity of $\vartheta[m]$ along $q_\nu=0$ in the coordinates
$$q_0=e^{2\pii (z_0+z_1)/8},\quad q_2=e^{2\pii (z_2+z_1)/8},
\quad q_1=e^{-2\pii z_1/8}.$$
is
$$a_1,\quad a_2,\quad a_1+a_2-2a_1a_2,\qquad m={a\choose b}.$$
}
Multt%
\finishproclaim
Let $m_1,\dots,m_4$ be a syzygetic quadruple of even characteristics
and let $\gotn=\{n_1,\dots,n_6\}$ be the complementary even characteristics. We
introduced the modular form
$$T=T_{\Gotn}=\vartheta[n_1]\cdots\vartheta[n_6].$$
By \Thst\ this is a cusp form of weight three for a conjugate of the 
group $\Gamma_0[2]$. Its
character is trivial on $\Gamma_2[4]$. 
Hence it is better now to use $l=4$, i.e.\ to change the notation
and to use the normal coordinates
$$q_0=e^{2\pii (z_0+z_1)/4},\quad q_2=e^{2\pii (z_2+z_1)/4},
\quad q_1=e^{-2\pii z_1/4}.$$
We can consider the differential form
$$\omega =Tdz_0\wedge dz_1\wedge dz_2$$
on $\tilde X(4)$.
For general reasons this gives a regular differential form on
the whole $\tilde X(4)$. We are interested in the zero divisor of this form.
We pull it back to $\tilde D$ to obtain
$$C{T\over q_0q_1q_2}dq_0dq_1dq_2.$$
At the moment we only are interested in the behavior of the differential forms
along the three divisors $q_\nu=0$.
From \Multt\ we can read off the vanishing order of $T$ along $q_\nu=0$.
Of course the result depends on the choice of the syzygetic
quadruple.
For example 
the vanishing order of the power series along $q_0=0$ is given by
$${a_{11}+\cdots+a_{61}\over 2},\qquad n_\nu={a_\nu\choose b_\nu}.$$
(The denominator $2$ occurs since we used  different $q_\nu$ in \Multt.)
By inspection of the 15 cases one sees that this expression is always 1 or 2.
Since for the order of differential form one has to subtract 1,
one gets for the vanishing
order of the differential form along $q_0$ either $0$ or $1$.
The same argument works for the variables $q_1$ and $q_2$.
Now one has to use a list of all the 15  syzygetic quadruples.
By inspection one finds:
\proclaim
{Proposition}
{There are\/ $15$ syzygetic quadruples.
We denote the vanishing order of the differential form
$\omega$ (pulled back to $\tilde D$)  at $q_\nu=0$ by $k_\nu$.
Then one has for $(k_0,k_1,k_2)$
the possibilities
\smallni
$(0,0,0)$ (eight cases),\hfill\break
$(1,1,1)$ (one case),\hfill\break
$(0,0,1)$ (two cases),\hfill\break
$(0,1,0)$ (two cases),\hfill\break
$(1,0,0)$ (two cases).\hfill\break
}
vOrd%
\finishproclaim
Let us assume that the order along $q_0=0$ is one. Then
we are in the case $a_{11}+\cdots+a_{61}=4$.
A glance at the power series shows that then  only even $t_0$ 
occur. Hence the series is 
invariant under the transformation $q_0\mapsto -q_0$.
This transformation is induced by the translation $z_0\mapsto z_0+2$.
Hence this translation belongs to our group $\Gamma_{\Gotn}$. This translation acts on
$\tilde X(4)$ as a reflection, meaning  that it is of order 2 and
fixes an irreducible subvariety of codimension 1.
So we obtain the following result which extends \NullD\ to the boundary:
\proclaim
{Proposition}
{The modular form $T$ induces on $\tilde X(4)$ a differential form with the following
property:
if $F\subset\tilde X$ is an irreducible component of its zero divisor
then there exists a reflection in $\Gamma/\pm\Gamma_2[4]$, with fixed point
set $F$. The multiplicities of the zeros are one.}
RaD%
\finishproclaim
We need the following general result:
\proclaim
{Theorem}
{Let $X$ be a quasi-projective smooth variety of dimension
three and $G$ a finite group of automorphisms
of it. Assume that every point of $M/G$
admits an open neighborhood such that
on its regular locus  exists a three-form without zeros.
Then $X/G$ admits a
quasi-projective crepant%
\footnote{*)}{\ninepoint\rm ``Crepant'' means that the inverse images of these
three-forms have no poles or zeros at the exceptional divisors.}
desingularization.}
Ito%
\finishproclaim
{\it Proof.\/} This theorem is a consequence of theorem 1.2 in [BKS].
Unfortunately the assumptions  are slightly different from ours.
This theorem only directly applies if an arbitrary stabilizer $G_a$ acts
on the tangent space as a subgroup of $\SL()$.
\smallskip
For this reason we want to give some details.
A cluster is a
zero dimensional closed subscheme of  $M$. The set of all
clusters of a fixed degree
is parametrized by a quasi-projective Hilbert scheme.
The set of all $G$-invariant clusters of degree $\#G$
is a closed subscheme. We denote by $G$-$\Hilb(M)$ the irreducible
component which contains the generic clusters (given by
$\#G$ distinct points.)
Then
$$G\hbox{-}\Hilb (M)\lo M$$
is a crepant resolution. This can be proved
analytically locally. Therefore we can assume
that $M=\cz^3$ and $G\subset\GL(n,\cz)$ is a linear
group. Consider the subgroup $G'\in G$ 
generated by all pseudo-reflections. By a result of
Chevalley $\cz^3/G'$ can be identified with $\cz^3$ in such a way
that $H=G/G'$ acts linearly. Now our assumption gives that
$H$ acts by a subgroup of $\SL()$. Hence we can apply the above
mentioned theorem 1.2 in [BKR] and obtain that
$$H\hbox{-}\Hilb (\cz^3/G')\lo \cz^3/G$$
is a crepant resolution.
Taking inverse images one obtains a natural closed embedding
$$H\hbox{-}\Hilb (\cz^3/G')\lo  G\hbox{-}\Hilb (\cz^3/G).$$
Since it is birational and since $G\hbox{-}\Hilb (\cz^3/G)$
is irreducible, it must be an isomorphism.\qed
\smallskip
In our main example,
concerning the subgroup of index two of $\Gamma_2[2]$,
we can avoid using
this rather deep theory, since the occurring singularities
are rather mild and easy
to desingularize. We will describe this below.
\smallskip
Our main result is:
\proclaim
{Theorem}
{Let be $\Gamma'$ be any group between $\Gamma_{\Gotn}$ and $\Gamma_0[2]_{\Gotn}$.
Then the Satake compactification of $\hz_2/\Gamma'$ admits a desingularization
which is a Calabi-Yau manifold.}
Main%
\finishproclaim
It remains to notice that the first Betti number vanishes.
It is known that $h^{1,0}$ is zero
for all non-singular models of Siegel modular varieties of genus $>1$.
From Serre duality, it follows that also $h^{2,0}=0$.
This fits to the computation in [LW] .
For $\Gamma_2[4]$ one has
$h^{2,0}=6$. It is easy to derive from the
description given in [LW] that any form  invariant under $\Gamma'$, is zero.\qed
\smallskip
We will give now an explicit and very simple construction 
of the desingularization in the case
of the smallest group $\Gamma=\Gamma_2[2]_{\Gotn}$.
We have to study the singularities of $\tilde X(4)/\Gamma$.
For this we introduce the finite groups
$$A=\Gamma_2[2]/\pm\Gamma_2[4]\quad\hbox{and}\quad  B=\Gamma/\pm\Gamma_2[4].$$
So $B$ is a subgroup of index two of $A$. The basic point now is that $A$ is an abelian
group and that each element of $A$ has order $\le 2$. Igusa [Ig2] proved that
$\tilde X(4)/A$ is smooth. Let $a\in \tilde X(4)$ be some point. Using suitable
coordinates the stabilizer $A_a$ can be be linearized. Hence $A_a$ can be considered
as some subgroup of $\GL(3,\cz)$. Since $A_a$ is abelian we can diagonalize it.
Hence $A_a$ can be considered as subgroup of $(\gz/2\gz)^3$, where this group acts
on $\cz^3$ in the obvious way by changing signs. Since the quotient by $A_a$ is smooth,
$A_a$ must be generated by reflections (here simply sign changes of one variable).
Since $B$ is a subgroup of index two of $A$, we know that $B_a$ is a subgroup
of index $\le 2$ of $A_a$. Hence we obtain:
\proclaim
{Lemma}
{Let $\Gamma=\Gamma_2[2]_{\Gotn}$. The quotient $\tilde X(4)/\Gamma$ looks locally
like $\cz^3/H$, where $H\subset(\gz/2\gz)^3$ is a subgroup which is contained
in some reflection subgroup as subgroup of index $\le 2$.}
LRs%
\finishproclaim
The reflection subgroups of $(\gz/2\gz)^3$ are trivial to describe.
Up to permutation of the variables one has four cases:
the trivial subgroup, sign change of the first variable, arbitrary sign changes of the
first two variables and sign changes of all variables.
\smallskip
The group $H$ is a subgroup of index $\le 2$ of such a group. We are only interested
in cases where $H$ is not generated by reflections.
There are only three types of such groups:
\smallni
\item{1)} The group of order two which is generated by
$$(z_1,z_2,z_3)\loma (-z_1,-z_2,z_3).$$
\item{2)} The group of order 4 which is generated by
$$(z_1,z_2,z_3)\loma (-z_1,-z_2,z_3)\quad\hbox{and}\quad (z_1,z_2,z_3)\loma (z_1,z_2,-z_3).$$
\item{3)} The group of order 4 given by sign changes of an even number
of coordinates.
\smallni
In the second case the group $H$ contains a reflection, changing the sign of
the third variable. We can take the quotient by this reflection and thus reduce
the second case to the first one.
So there only two types of singularities occur, case 1) and 3). In both cases the
differential form can be written as $h(z)dz_1\wedge dz_2\wedge dz_3$
with $h$ some $H$-invariant holomorphic function in a neighborhood of the origin.
Since this differential form should produce a differential form without zeros
on the regular locus of the quotient, we obtain that $h(0)\ne 0$.
This means that we can assume that $\omega$ simply is given by
$dz_1\wedge dz_2\wedge dz_3$.
\smallskip 
Now we describe the desingularization.
In  case 1) we blow up $\cz^3$
along the line $z_1=z_2=0$. A typical affine chart of the blow-up is given by
$(w_1=z_1/z_2,z_2,z_3)$. Because of
$$z_2dz_1\wedge dz_2\wedge dz_3=dw_1\wedge dz_2\wedge dz_3$$
the pull-back of the differential form $dz_1\wedge dz_2\wedge dz_3$
has a zero along $z_2=0$.
The group $H$ acts now as the group generated by $(w_1,z_2,z_3)\mapsto (w_1,-z_2,z_3)$.
This is a reflection group. Hence the quotient is smooth and the zero
of $\omega$ along the ramification $z_2=0$ disappears on the quotient, which gives the desired
desingularization of $\cz^3/H$.
\smallskip
The remaining case 3) is slightly more involved.
The group $H$ now has order 4 and consists of arbitrary sign
changes of an even number of variables.
The singular locus is the image of the union of the three coordinate axes in $\cz^3$.
This time there is no canonical way to desingularize! We have to make a choice.
We choose one of the three coordinate axes, for example $z_1=z_2=0$. 
We start by blowing
up this line.
A typical affine
chart of the blow up again
is $(w_1=z_1/z_2,z_2,z_3)$. The differential form gets a zero along $z_2=0$.
Now we have to consider the strict transform of the singular locus. It is given
by $w_1=z_3=0$ We have to blow up this locus up again. A typical chart is
$(u_1=w_1/z_3,z_2,z_3)$. The differential form $dz_1\wedge dz_2\wedge dz_3$
now gets besides $z_2=0$ the additional zero
$z_3=0$. The group acts on the coordinates $(u_1=z_1/(z_2z_3),z_2,z_3)$ by arbitrary
sign changes of the variables $z_2,z_3$. This is a reflection group. The quotient is smooth
and the zeros of the differential form disappear on the quotient.
\smallskip
Since in the third case we have no canonical resolution
(different choices of the coordinate axes lead to resolutions which are related by
flops), we have to explain how to glue the resolutions to get a resolution
of the global $\tilde X(4)/\Gamma$. The point is that the singular locus of $\tilde X(4)/\Gamma$
contains itself only finitely many singular points. These points lead to the case 3).
For each of these finitely many points one has to make a choice. But  the
smooth points of the singular locus lead to case 1) where we have a canonical resolution.
Therefore everything fits together to a complex manifold.
(This argument does not give projectivity.)\qed
\neupara{Equations}%
In this section we treat the case where $\gotn$ is complementary to the standard
syzygetic quadruple and $\Gamma$ means $\Gamma_{\Gotn}$ for this choice.
\smallskip
We give the equations for the variety $X$ in the Igusa coordinates.
To simplify the notation we will write the theta constants in the form
$$\vartheta[m]=\vartheta\Bigl[{a_1\;a_2\atop b_1\;b_2}\Bigr]\quad
\hbox{for}\quad m=\pmatrix{a_1\cr a_2\cr b_1\cr b_2}.$$
The results which we describe now can by taken from Igusa's paper [Ig1]
(see page 397).
\proclaim
{Theorem (Igusa)}
{The five modular forms
$$y_0=\vartheta\Bigl[{ 00\atop 11}\Bigr]^4,\quad
y_1=\vartheta\Bigl[{00\atop01}\Bigr]^4,\quad
y_2=\vartheta\Bigl[{00\atop00}\Bigr]^4,$$
$$y_3=-\vartheta\Bigl[{10\atop00}\Bigr]^4- \vartheta\Bigl[{00\atop11}\Bigr]^4,\quad
y_4= -\vartheta\Bigl[{10\atop01}\Bigr]^4- \vartheta\Bigl[{00\atop11}\Bigr]^4$$
generate the ring of modular forms of even weight (with trivial multipliers)
for $\Gamma_2[2]$. The defining relation is the
the quartic equation
$$(y_0y_1+y_0y_2+ y_1y_2- y_3y_4)^2= 4y_0y_1y_2( y_0+y_1+y_2+
y_3+y_4).$$}
Igz%
\finishproclaim
Moreover according to Igusa we have the relation
$$ 2\vartheta\Bigl[{00\atop01}\Bigr]^2 \vartheta\Bigl[{00\atop00}\Bigr]^2
\vartheta\Bigl[{00\atop10}\Bigr]^2 \vartheta\Bigl[{00\atop11}\Bigr]^2
= y_0y_1+y_0y_2+ y_1y_2- y_3y_4.$$
We set
$$y_5:=\Theta\qquad 
\bigl(= \vartheta\Bigl[{00\atop01}\Bigr] \vartheta\Bigl[{00\atop00}\Bigr]
\vartheta\Bigl[{00\atop10}\Bigr] \vartheta\Bigl[{00\atop11}\Bigr]\bigr).$$
We recall that the Calabi-Yau form $T$ is the product of the complementary six thetas.
The product $y_5T$ is Igusa's modular form $\chi_5$ which is the unique modular form
of weight 5 for the full modular group. The character of this form is non-trivial but
trivial on $\Gamma_2[2]$. Hence $y_5$ is a modular form of weight two with
trivial character on our group $\Gamma$.
\proclaim
{Proposition}
{The ring of modular forms of even weight for $\Gamma$ is generated
by the six forms of weight two, $y_0,\dots,y_5$. The defining
relations are the quartic
$$y_5^4=y_0y_1y_2( y_0+y_1+y_2+ y_3+y_4)$$
and the quadric
$$2y_5^2= y_0y_1+y_0y_2+ y_1y_2- y_3y_4.$$
}
RiGaz%
\finishproclaim
{\it Proof.\/}
The field of modular functions for $\Gamma$ is a quadratic extension
of the field of modular functions for $\Gamma_2[2]$. Hence the homogenous
field of fractions of $\cz[y_0,\dots,y_5]$ is the full field of modular functions
of $\Gamma$. Since $y_0,\dots,y_5$ have no common zero on the Satake compactification,
the ring of all modular forms of even weight is the normalization of 
$\cz[y_1,\dots,y_6]$. Hence it suffices to show that the ideal given by the above
two relations is a prime ideal and that the quotient is a normal ring.
Since we have two relations, the factor ring is a complete intersection
and hence a Cohen-Macaulay ring.
Since the singular locus (as has been described in the introduction) is of codimension
$\ge 2$, on can apply the well-known Serre criterion for normality.\qed
\neupara{The Calabi-Yau form}%
Now we use Runge's approach and consider
the theta series of second kind
$$f_a(Z):=\vartheta\Bigl[{a\atop 0}\Bigr](2Z).$$
They are linked to the $\vartheta[m]$ by the classical relation
$$ \vartheta\Bigl[{a\atop b}\Bigl](Z)^2=\sum_{x\in\{0,1\}^2}(-1)^{b'x}f_{a+x}(Z)f_x(Z).$$
Hence the rings $\cz[f_af_b]$ and $\cz[\vartheta[m]^2]$ agree.
We denote the $f_a$  by $f_1,f_2,f_3,f_4$ in the ordering
$$\pmatrix{0\cr 0}\quad\pmatrix{1\cr 0}\quad\pmatrix{0\cr 1}\quad\pmatrix{1\cr 1}.$$
The following modified version of Igusa's result \Igz\ is due to Runge [Ru]:
\proclaim
{Proposition}
{The algebra of modular forms of even weight (with trivial multipliers) with respect
to the group $\Gamma_2[2]$ is generated by the five forms, all of weight two,
$$\eqalign{
F_1&= f_1^4 + f_2^4 + f_3^4 + f_4^4,\cr
F_2&= f_1^2f_2^2 + f_3^2f_4^2,\cr
F_3&= f_1^2f_3^2 + f_2^2f_4^2,\cr
F_4&= f_1^2f_4^2 + f_2^2f_3^2,\cr
F_5&= f_1f_2f_3f_4.\cr}
$$
The defining relation is
$$16F_5^4=
-F_1^2F_5^2+F_1F_2F_3F_4-F_2^2F_3^2-F_2^2F_4^2+4F_2^2F_5^2-F_3^2F_4^2+
4F_3^2F_5^2+4F_4^2F_5^2.$$
}
tilAG%
\finishproclaim
It is very easy to describe the action of $\Gamma_{2,0}[2]$ on the generators:
\proclaim
{Lemma}
{The three translations
$$Z\loma Z+S;\quad S= \pmatrix{1&0\cr0&0},\
\pmatrix{0&1\cr1&0},\ \pmatrix{0&0\cr0&1}$$
act by $(F_1,\dots,F_5)\mapsto$
$$ (F_1,-F_2,F_3,-F_4,F_5),\quad
(F_1,F_2,F_3,F_4,-F_5),\quad(F_1,F_2,-F_3,-F_4,F_5).$$
The unimodular substitutions
$$Z\loma Z[U];\quad U=\pmatrix{0&1\cr 1&0},\quad U=\pmatrix{1&1\cr 0&1}$$
act by
$$(F_1,\dots,F_5)\loma (F_1,F_3,F_2,F_4,F_5),\quad
(F_1,F_2,F_4,F_3,F_5).$$}
ActGn%
\finishproclaim
The analogue of \RiGaz\ is:
\proclaim
{Proposition}
{The ring of modular forms of even weight (with trivial multipliers) of $\Gamma$
is generated by six modular forms all of weight two, namely $F_1,\dots,F_5$ (see \tilAG)
and the additional form $F_6=\Theta$. There are two defining relations, namely
the relation described in \tilAG\ and the additional quadratic relation
$$F_6^2=F_1^2-4F_2^2-4F_3^2-4F_4^2+32F_5^2.$$
}
RiGa%
\finishproclaim
We know  the action of $\Gamma_{2,0}[2]$ also on $F_6$ (\LeRi).
We give just one example:
\proclaim
{Lemma}
{The matrix
$$M=\pmatrix{E&0\cr C&E}\quad C=\pmatrix{2&2\cr 2&2}$$
acts by
$$(F_1,\dots,F_6)\loma (F_1,F_2,F_3,F_4,F_5,-F_6).$$}
Trafsp%
\finishproclaim
In the new coordinates we compute an algebraic expression for the Calabi-Yau differential form
$$\omega:=T dz_0\wedge dz_1\wedge dz_2,\qquad Z=\pmatrix{z_0&z_1\cr z_1&z_2}.$$
As we mentioned
$Ty_5=\chi_5$ is the well-known cusp form of weight
5 for the full modular group.
We will use the homogenous Jacobian
$$W(f_1,f_2,f_3,f_4)=
\det\pmatrix{f_1&f_2&f_3&f_4\cr
\partial_0f_1&\partial_0f_2&\partial_0f_3&\partial_0f_4\cr
\partial_1f_1&\partial_1f_2&\partial_1f_3&\partial_1f_4\cr
\partial_2f_1&\partial_2f_2&\partial_2f_3&\partial_2f_4\cr}.$$
Here $\partial_i$ denotes differentiation by $z_i$.
The connection with the usual Jacobian
$$J(f_1/f_4,f_2/f_4,f_3/f_4)=
\det\pmatrix{
\partial_0(f_1/f_4)&\partial_0(f_2/f_4)&\partial_0(f_3/f_4)\cr
\partial_1(f_1/f_4)&\partial_1(f_2/f_4)&\partial_1(f_3/f_4)\cr
\partial_2(f_1/f_4)&\partial_2(f_2/f_4)&\partial_2(f_3/f_4)\cr
}$$
is
$$W(f_1,f_2,f_3,f_4)=f_4^4J(f_1/f_4,f_2/f_4,f_3/f_4).$$
The Jacobian of a modular substitution $M$ is $\det(CZ+D)^{-3}$.
Hence $J$ is a modular form of weight $3$ and $W$ is a modular form of weight
5. If one applies a  modular transformation to the $f_i$ one obtains a linear
transformation of them. This shows that $W$ is invariant under $M$ up to the determinant
of this linear transformation. This shows that $W$ up to a constant factor equals
Igusa's modular form
$\chi_5$ which is the only modular
form of weight 5 for the full modular group.
It can be defined as the product of the ten theta series.
As a consequence we obtain
$$d(f_1/f_4)\wedge d(f_2/f_4)\wedge d(f_3/f_4)=
c\;{\chi_5\over f_4^4}\; dz_0\wedge dz_1\wedge dz_2$$
with a certain constant $c$. Using $T\Theta=\chi_5$ we get
$$c\omega={f_4^4\over\Theta}\;d(f_1/f_4)\wedge d(f_2/f_4)\wedge d(f_3/f_4).$$
We set
$$g_1={f_1\over f_4},\quad g_2={f_2\over f_4},\quad g_3={f_3\over f_4}$$
and, with the notations of \tilAG
$$G_1={F_2\over F_5},\quad G_2={F_3\over F_5},\quad G_3={F_4\over F_5}$$
We get
$$G_1={g_1g_2\over g_3}+{g_3\over g_1g_2},\quad
G_2={g_1g_3\over g_2}+{g_2\over g_1g_3},\quad
G_3={g_2g_3\over g_1}+{g_1\over g_2g_3}.$$
The Jacobian of this rational transformation is
$$4{(g_3^2-g_1^2g_2^2)(g_2^2-g_1^2g_3^2)(g_1^2-g_2^2g_3^2)\over g_1^4g_2^4g_3^4}$$
or
$$4{(f_1^2f_2^2-f_3^2f_4^2)(f_1^2f_3^2-f_2^2f_4^2)(f_1^2f_4^2-f_2^2f_3^2)
\over f_1^4f_2^4f_3^4}$$
This gives
$$\eqalign{
4c\omega
&={f_1^4f_2^4f_3^4f_4^4\;d(F_2/F_5)\wedge d(F_3/F_5)\wedge d(F_4/F_5)
\over\Theta (f_1^2f_2^2-f_3^2f_4^2)
(f_1^2f_3^2-f_2^2f_4^2)(f_1^2f_4^2-f_2^2f_3^2)}
\cr
&={F_5^4\over (F_2F_3F_4-2F_1F_5^2)F_6}\;
d(F_2/F_5)\wedge d(F_3/F_5)\wedge d(F_4/F_5).\cr}$$
The group $\Gamma_{\Gotn}$ is normal in $\Gamma_{2,0}[2]$. The quotient acts
on the ring of modular forms for $\Gamma_{\Gotn}$. The subgroup of index two
$$G=\Gamma_{2,0}[2]_{\Gotn}/\Gamma_2[2]_{\Gotn}$$
leaves the differential form invariant. Using \Trafsp\ and \ActGn\ we can express
the action of $G$ on the generators. We obtain:
\proclaim
{Theorem}
{Let $X$ be the subvariety of $P^5(\cz)$ given a as intersection
of the quartic
$$16x_4^4+x_0^2x_4^2+x_1^2x_2^2+x_1^2x_3^2+x_2^2x_3^2=x_0x_1x_2x_3+4x_4^2(x_1^2+x_2^2+x_3^2)$$
and the quadric
$$x_5^2=x_0^2-4x_1^2-4x_2^2-4x_3^2+32x_4^2.$$
This is a normal projective variety
of dimension three. There exists a desingularization
$\tilde X\to X$ which is a Calabi-Yau manifold.
The differential form without zeros is given by
$$\omega:={x_4^4\over (x_1x_2x_3-2x_0x_4^2)x_5}\;
d(x_1/x_4)\wedge d(x_2/x_4)\wedge d(x_3/x_4).$$
The group $G\cong S_3\cdot (\gz/2\gz)^3$ (semidirect product) generated by
\vskip1mm\item{\rm 1)} arbitrary
permutations of $x_1,x_2,x_3$ followed by the sign change of $x_5$ if
the permutation is odd,
\item{\rm 2)} arbitrary sign changes of two of the $x_1,x_2,x_3$,
\item{\rm 3)}
the sign change of $x_4$,
\smallni
is a group of automorphisms of $X$ which fixes $\omega$.
For each subgroup $H\subset G$ the quotient $X/H$ admits a Calabi-Yau desingularization
with Calabi-Yau form $\omega$.
}
VaR%
\finishproclaim
For sake of completeness we give the coordinate transformation between the
coordinates $x_i$ and the coordinates $y_i$ which we used in the introduction.
We have $y_5=x_5$ and
$$\pmatrix{y_0\cr y_1\cr y_2\cr y_3\cr y_4}=
\pmatrix{
1&-2&-2&2&0\cr1&-2&2&-2&0\cr1&2&2&2&0\cr-1&2&-2&-2&-8\cr-1&2&-2&-2&8\cr}
\pmatrix{x_0\cr x_1\cr x_2\cr x_3\cr x_4}$$
The expression we obtained for $\omega$ and for the action of
$G$ in the coordinates $y_i$ looked not
very nice, hence we skip them.
\vskip1cm\noindent
{\paragratit References}%
\bigni
\item{[BN]} Barth, W. Nieto, I.:  {\it Abelian surfaces of type $(1,3)$ and
quartic surfaces with 16 skew lines,\/}
J. Algebraic Geometry {\bf 3}, 173--222 (1994)
\medskip
\item{[BKR]}
Bridgeland, T., King, A., Reid, M.:
{\it The McKay correspondence as an equivalence of derived categories,\/}
J.~Amer.~Math.~Soc. {\bf 14}, 535–-554 (2001)
\medskip
\item{[Fr]} Freitag, E.: {\it Siegelsche Modulfunktionen,} Grundlehren
der mathematischen Wissenschaften, Bd. {\bf 254}. Berlin Heidelberg New
York: Springer (1983)
\medskip
\item{[GC]} Gritsenko, V., Cl\'ery, F.:
{\it The Siegel modular form of genus 2 with the simplest divisor,\/}
arXiv: 0812.3962 (2008)
\medskip
\item{[GH]} Gritsenko, V., Hulek, K.:
{\it The modular form of the Barth-Nieto quintic.\/}
IMRN International Math. Research Notices, No {\bf 17} (1999)
\medskip
\item{[GN]} van Geemen, B., Nygaard, N.O.:
{\it On the geometry and arithmetic of some Siegel modular threefolds,\/}
Journal of Number Theory {\bf 53}, 45--87 (1995)
\medskip
\item{[GS]} van Geemen, B., van Straten, D.: {\it The cusp form of weight 3 on
$\Gamma_2(2,4,8)$,\/}
Mathematics of Computation {\bf 61}, 849--872 (1993)
\medskip
\item{[GHSS]} van Geemen, B., Hulek, K., Spandau, J., van Straten, D.:
{\it The modularity of the Barth-Nieto quintic and its relatives,\/}
Advances in Geometry {\bf 1}, 263--289 (2001)
\medskip
\item{[Ig1]} Igusa, I.: {\it On Siegel modular forms of genus two II,\/}
Amer. J. Math. {\bf 86}, 392--412 (1964)
\medskip
\item{[Ig2]} Igusa, I.: {\it A desingularization problem in the theory of
Siegel modular functions,\/}  Math. Annalen {\bf 168} 228--260
(1967)
\medskip
\item{[Ito]} Ito, Y.:
{\it Gorenstein quotient singularities of monomial type in dimension three,\/}
J. Math. Sci. Univ. Tokyo {\bf 2}, 419--440 (1995)
\medskip
\item{[LW]} Lee, R., Weintraub, S.H.:
{\it The Siegel modular variety of degree two and level four,\/}
Mem. Am. Math. Soc. {\bf 631}, 3-58 (1998)
\medskip
\item{[Re]} Reid, M.: {\it La correspondence de McKay,\/}
S\'eminaire Bourbaki 1999/2000. Ast\'erisque No. {\bf 276}, 53--72 (2002)
\medskip
\item{[Ru1]} Runge, B.:
{\it On Siegel modular form, part I\/}, J. Reine angew. Math. {\bf 436}, 57-85 (1993)
\medskip
\item{[Ru2]} Runge, B.:
{\it On Siegel modular forms, part II\/}, Nagoya Math. J. {\bf138}, 179-197 (1995)
\medskip
\item{[SM]} Salvati Manni, R.:
{\it Thetanullwerte and stable modular forms,\/}
Am. J. of Math. {\bf 111}, 435--455 (1989)
\bye